\documentclass[a4paper,12pt]{amsart}
\usepackage{amssymb}
\usepackage{ifthen}
\usepackage{graphicx}
\usepackage{float}
\usepackage{caption}
\usepackage{subcaption}
\usepackage{cite}
\usepackage{amsfonts}
\usepackage{amscd}
\usepackage{amsxtra}
\usepackage{color}

\newtheorem{thm}{Theorem}
\newtheorem{cor}{Corollary}
\newtheorem{lem}{Lemma}
\newtheorem{rem}{Remark}

\newtheorem{conj}{Conjecture}
\newtheorem{prob}{Problem}

\theoremstyle{definition}
\newtheorem{defn}{Definition}[section]
\newtheorem{example}{Example}

\newenvironment{pf}[1][]{%
 \vskip 1mm
 \noindent
 \ifthenelse{\equal{#1}{}}%
  {{\slshape Proof. }}%
  {{\slshape #1.} }%
 }%
{\qed\bigskip}

\newcounter{alphabet}
\newcounter{tmp}
\newenvironment{Thm}[1][]{\refstepcounter{alphabet}%
\bigskip%
\noindent%
{\bf Theorem \Alph{alphabet}}%
\ifthenelse{\equal{#1}{}}{}{ (#1)}%
{\bf .} \itshape}{\vskip 8pt}

\makeatletter
\newcommand{\Ref}[1]{\@ifundefined{r@#1}{}{\setcounter{tmp}{\ref{#1}}\Alph{tmp}}}
\makeatother

\newcommand{\IN}{{\mathbb N}}

\newcommand{\ID}{{\mathbb D}}

\def\be{\begin{equation}}
\def\ee{\end{equation}}

\newcommand{\bee}{\begin{enumerate}}
\newcommand{\eee}{\end{enumerate}}

\newcommand{\blem}{\begin{lem}}
\newcommand{\elem}{\end{lem}}
\newcommand{\bthm}{\begin{thm}}
\newcommand{\ethm}{\end{thm}}
\newcommand{\bcor}{\begin{cor}}
\newcommand{\ecor}{\end{cor}}
\newcommand{\beg}{\begin{example}}
\newcommand{\eeg}{\end{example}}
\newcommand{\begs}{\begin{examples}}
\newcommand{\eegs}{\end{examples}}
\newcommand{\bdefe}{\begin{defn}}
\newcommand{\edefe}{\end{defn}}
\newcommand{\bprob}{\begin{prob}}
\newcommand{\eprob}{\end{prob}}
\newcommand{\bques}{\begin{ques}}
\newcommand{\eques}{\end{ques}}
\newcommand{\bei}{\begin{itemize}}
\newcommand{\eei}{\end{itemize}}
\newcommand{\bcon}{\begin{conj}}
\newcommand{\econ}{\end{conj}}
\newcommand{\bcons}{\begin{conjs}}
\newcommand{\econs}{\end{conjs}}
\newcommand{\bprop}{\begin{propo}}
\newcommand{\eprop}{\end{propo}}
\newcommand{\br}{\begin{rem}}
\newcommand{\er}{\end{rem}}
\newcommand{\brs}{\begin{rems}}
\newcommand{\ers}{\end{rems}}
\newcommand{\bo}{\begin{obser}}
\newcommand{\eo}{\end{obser}}
\newcommand{\bos}{\begin{obsers}}
\newcommand{\eos}{\end{obsers}}
\newcommand{\bpf}{\begin{pf}}
\newcommand{\epf}{\end{pf}}
\newcommand{\ba}{\begin{array}}
\newcommand{\ea}{\end{array}}
\newcommand{\beq}{\begin{eqnarray}}
\newcommand{\beqq}{\begin{eqnarray*}}
\newcommand{\eeq}{\end{eqnarray}}
\newcommand{\eeqq}{\end{eqnarray*}}

\newcommand{\Llra}{\Longleftrightarrow}

\newcommand{\ds}{\displaystyle}

\newcounter{minutes}
\divide\time by 60
\newcounter{hours}
\multiply\time by 60 \addtocounter{minutes}{-\time}
\begin{document}

\bibliographystyle{amsplain}

%

\title[Bohr type inequalities for functions with a multiple zero at the origin] {Bohr type inequalities for functions with a multiple zero at the origin   }

\thanks{
File:~\jobname .tex,
          printed: \number\day-\number\month-\number\year,
          \thehours.\ifnum\theminutes<10{0}\fi\theminutes}

\author[S. Ponnusamy]{Saminathan Ponnusamy
}
\address{
S. Ponnusamy, Department of Mathematics,
Indian Institute of Technology Madras, Chennai-600 036, India.
}
\email{samy@iitm.ac.in}
%

\author[K.-J. Wirths]{Karl-Joachim Wirths}
\address{K.-J. Wirths, Institut f\"ur Analysis und Algebra, TU Braunschweig,
38106 Braunschweig, Germany.}
\email{kjwirths@tu-bs.de}


\begin{abstract}
Recently, there has been a number of good deal of research on the Bohr's phenomenon in various setting including a refined formulation of
his classical version of the inequality. Among them, in \cite{PaulPopeSingh-02-10} the authors considered the cases in which the above functions have a multiple zero at the origin. In this article, we present a refined version of Bohr's inequality for these cases and give a partial answer to a question from \cite{PaulPopeSingh-02-10} for the revised setting.\\

{\bf (Dedicated to the memory of Professor Stephan Ruscheweyh)}
\end{abstract}

\maketitle
\pagestyle{myheadings}
\markboth{S. Ponnusamy and K.-J. Wirths}{ Bohr Inequality in case of a multiple zero}

\section{Introduction and Statement of Results}

Let  $\mathcal B$ denote the class of all analytic functions $f$  defined on  the open unit disk $\mathbb{D}:=\{z\in \mathbb{C}:\,|z|<1\}$
such that $|f(z)|\leq 1$ in $\ID$.   The classical inequality of H. Bohr \cite{Bohr-14} in 1914
asserts that if  $f\in {\mathcal B}$ and $f(z)=\sum_{n=0}^{\infty} a_n z^n,$
then
\be\label{KKP1-eq1}
\sum_{n=1}^{\infty} |a_n| r^n \leq 1 -|a_0| 
~\mbox{ for  $r\leq 1/3$}
\ee
and the number $1/3$, which is called the Bohr radius for the family ${\mathcal B}$, cannot be improved. We would
like to point out that Bohr originally established the inequality \eqref{KKP1-eq1} only for $r\leq 1/6$ and the value $1/3$
was obtained independently by M. Riesz, I. Schur and N. Wiener.
For more information about Bohr's inequality and related investigations, we refer
to the recent expository articles  \cite{AAP2016} and \cite[Chapter 8]{GarMasRoss-2018}.
There are many proofs of this inequality (cf. \cite{Sidon-27-15} and \cite{Tomic-62-16}). Indeed, we may assume
without loss of generality that $|f(z)|< 1$ in $\ID$, i.e.
$$f(z)\prec \varphi_{a_0} (z), \quad  \varphi_ a(z)=\frac{a-z}{1-\overline{ a}z }=a -(1-|a|^2)\frac{z}{1-\overline{a}z },
$$
where $\prec$ is the usual subordination, and  for each $a \in \ID$, $\varphi_a (\ID)=\ID$ is a convex domain. It follows from \cite{Rogo-43} that
$|a_n| \leq |\varphi_{a_0} '(0)|=1-|a_0|^2$ for $n\geq 1$ and hence,
$$\sum_{n=1}^{\infty} |a_n| r^n \leq (1 -|a_0|^2)\frac{r}{1-r}\leq  (1 -|a_0|)\frac{2r}{1-r} \leq 1 -|a_0| ~\mbox{ for  $r\leq 1/3$}.
$$
The fact that $1/3$ cannot be improved follows from  $\varphi_a (z)$ because for this function, we have
$$\sum_{n=0}^{\infty} |a_n| r^n =1 + \frac{1-|a|}{1-|a| r}\left ( r(1+2|a|)-1\right ) >1  ~\Llra~ r >\frac{1}{1+2|a|},
$$
and, since $|a|$ can be chosen arbitrarily close to $1,$  the Bohr radius for $\mathcal{B}$ cannot be bigger than $1/3$. This simple result
attracted the attention of many to develop what is called Bohr-phenomenon in various setting including multi-dimensional analog of
it. See for example, \cite{AlKayPON-19,Bom-62, BDK5,BhowDas-18,BoasKhavin-97-4,BomBor-04,EvPoRa-19,GarMasRoss-2018,KayPon1, KayPon3, KayPon2,KP-AASFM2-19}.
Generalizations and extensions of Bohr's result can be found from  \cite{BDK5,PaulPopeSingh-02-10,PauSin--09} and the references therein.

Another natural generalization of the Bohr phenomenon is due to Fournier and Ruscheweyh \cite{FR}
who have considered the problem of determining the Bohr constant $B=B(D)$ such that
\be\label{bohr}
B=\sup \left \{r\in (0,1) : \sum_{k=0}^\infty|a_k|r^k\le 1
\mbox{ for all }~f(z):=\sum_{k=0}^\infty a_k r^k,\; z\in \mathbb{D} \right \},
\ee
where the supremum is taken over all functions $f$ analytic in $D$ with $f(D)\subset \overline{\mathbb{D}}$, and
$D\subset \mathbb{D}$ is a simply connected domain.

\begin{Thm} {\rm (Fournier and Ruscheweyh \cite{FR})} \label{FR-thm1}
Let $D_\gamma$ denote the disk $\{z\in \mathbb{C}:\,  \big|z+\frac{\gamma}{1-\gamma}\big|<\frac{1}{1-\gamma}\}$, for $0\leq \gamma<1$ and let $f:D_\gamma \to \mathbb{\overline{D}}$ be an analytic function such that $f(z)=\sum_{n=0}^\infty a_nz^n$ in $\mathbb{D}$. Then $B(D_\gamma)=\frac{1+\gamma}{3+\gamma}$ and $\sum_{n=0}^\infty |a_n|B(D_\gamma)^n=1$ if and only if $f\equiv c$, with $|c|=1$.
\end{Thm}

Recently, this result has been generalized to harmonic mappings in \cite{EvPoRa-20}.

It is worth pointing out that if $|a_0|$ in \eqref{KKP1-eq1} is replaced by $|a_0|^2$,
then the constant $1/3$ could be replaced by $1/2$. This fact also follows from examining the first inequality in the above inequalities, since $r/(1-r)\leq 1$ for
$r\in [0,1/2]$. Naturally, one can obviously expect a better estimate for functions which vanish at the origin up to some order. In view of this fact,
in \cite{PaulPopeSingh-02-10}, the authors considered among others for $k\in\IN$ the classes
\[{\mathcal B}_k:=\{f\in {\mathcal B}:\, f(0)= \cdots =f^{(k-1)}(0)=0\}
\]
and asked for which $r\in [0,1)$ and
\be\label{a}
f(z)=\sum_{n=k}^{\infty}a_nz^n\in {\mathcal B}_k
\ee
the inequality
\be\label{b}
\sum_{n=k}^{\infty}|a_n|r^n\leq 1
\ee
is valid.

For $k=1$,  Tomi\'c \cite{Tomic-62-16} proved that \eqref{b} holds for $0\leq r\leq 1/2$ (also obtained by Landau independently,
see \cite{Landau-29}). Later Ricci \cite{Ricci-1956} established that this holds for $0\leq r\leq 3/5$, and the largest value
of $r$ for which \eqref{b} holds would be in the interval $(3/5,1/\sqrt{2}]$. Later in 1962, Bombieri \cite{Bom-62}  found
that the inequality \eqref{b} holds for $r\in [0,1/\sqrt{2}]$, where the upper bound cannot be improved.
See \cite{KayPon1,KayPon2,KP-AASFM2-19} for new proofs of it in a general form. For $k\geq 2$, the authors
in \cite{PaulPopeSingh-02-10}  got the following partial answer.

\begin{Thm} \label{PVW1-theB}
{\rm (\cite[ Remark 2]{PaulPopeSingh-02-10})} Let $m(r):=\inf\left\{M(r), 1/\sqrt{1-r^2}\right\}$, where
\[    M(r)=\left \{ \begin{array}{rl} 1 & \mbox{ for $0\leq r \leq 1/3$}\\
\ds \frac{1-2r+5r^2}{4r(1-r)} & \mbox{ for $1/3 < r <1$}.
\end{array} \right .
\]
If $f\in \mathcal{B}_k$ has the expansion \eqref{a}, then \eqref{b} is valid for $r\in [0,r_k]$, where $r_k$ denotes the least positive root of the equation
$r^km(r)=1.$
\end{Thm}

In \cite{PaulPopeSingh-02-10}, the authors posed the problem to decide whether for $k\geq 2$ the upper bound $r_k$ is sharp.
For a direct comparison with main results of the paper, we refer to Table \ref{Table1} in which we indicate the numerical
values of $r_k$ for certain values of $k\geq 2$.

\begin{table}
\begin{center}
\begin{tabular}{|l|l||l|l||l|l||l|l||l|l||}
\hline
 {\bf $k$}   & $r_k$ & $k$ &  $r_k$  & $k$ & $r_k$ & {\bf $k$}   & $r_k$ & $k$ &  $r_k$\\
  \hline
2  & 0.786151 & 3 & 0.826031 & 4 &  0.851171 & 5  & 0.868837 & 6 & 0.882094   \\
  \hline
7  & 0.892493 & 8 & 0.900917 & 9 &  0.90791 & 10  & 0.913827 & 11 & 0.918911   \\
  \hline
15  & 0.933783 & 20 & 0.94546 & 25 & 0.953256 & 30  & 0.958885 & 35 & 0.963169   \\
  \hline
40  & 0.966553 & 50 & 0.97159 & 60 & 0.975183 & 70  & 0.977892 & 100 & 0.98315   \\
  \hline
\end{tabular}
\end{center}
\caption{$r_k$ is the unique root of the equation
$r^km(r)=1$  in $(0,1)$\label{Table1}}
\end{table}

The aim of the present paper is twofold. We want to give a number of sharp modified Bohr type inequalities for
functions in $\mathcal{B}_k$. With their help we will be able to give  answer to the above question in the
modified setting. The inequalities proved in the sequel are based on the following theorem.

\begin{Thm} {\rm(\cite{PVW})}\label{PVW1-theC} Let $f\in \mathcal{B}$ have the expansion $f(z)=\sum_{n=0}^{\infty} a_n z^n$, then
\[\sum_{n=0}^{\infty}|a_n|r^n+\left(\frac{1}{1+|a_0|}+\frac{r}{1-r}\right)\sum_{n=1}^{\infty}|a_n|^2r^{2n}\leq |a_0|+\frac{r}{1-r}(1-|a_0|^2).
\]
\end{Thm}

The proof of this theorem was given in \cite[Proof of Theorem 1]{PVW}. The main tool of this proof is a theorem due to Carlson, see \cite{Carlson-40} (and also \cite{PVW} for a simplified form of the proof). Moreover,
the case $f\in \mathcal{B}_1$ has been considered in \cite{PVW}, too. We use the same method to prove the following theorems for $f\in \mathcal{B}_k$, $k\geq 2.$

\bthm\label{PW1-th1}
For $k\geq 2$, let $f\in \mathcal{B}_k $ have an expansion  \eqref{a}. Then the inequality
\be\label{PW1-eq20-1}
\sum_{n=k}^{\infty}|a_n|r^n + \left(\frac{r^{-k}}{1+|a_k|}+\frac{r^{1-k}}{1-r}\right)\sum_{n=k+1}^{\infty}|a_n|^2r^{2n}\leq 1
\ee
is valid for $r\in [0,R_k]$, where $R_k$  is the unique root in $(0,1)$  of the equation
\[   4(1-r)-r^{k-1}(1-2r+5r^2)=0.\]
The upper bound $R_k$  cannot be improved.
\ethm

\br From the inequality \eqref{PW1-eq20-1}, it is clear that $R_k\leq r_k$.
\er

\bcor\label{PW1-cor1}
Let $k\geq 2$.  For any $f\in \mathcal{B}_k$ and any  $r\in [0,R_k]$, the strict inequality
\[\sum_{n=k}^{\infty}|a_n|r^n<1
\]
is valid.
\ecor

Numerical values of $R_k$ for certain values of $k\geq 2$ are presented in Table \ref{Table2}.
\begin{table}
\begin{center}
\begin{tabular}{|l|l||l|l||l|l||l|l||l|l||}
\hline
 $k$   & $R_k$ & $k$ &  $R_k$  & $k$ & $R_k$ & $k$   & $R_k$ & $k$ &  $R_k$\\
\hline
2  & 0.674837 & 3 & 0.720449 & 4 &  0.752379 & 5  & 0.776409 & 6 & 0.795346   \\
  \hline
7  & 0.81076 & 8 & 0.823614 & 9 &  0.834537 & 10  & 0.84396 & 11 & 0.852191   \\
  \hline
15  & 0.876981 & 20 & 0.897193 & 25 & 0.911051 & 30  & 0.921238 & 35 & 0.92909   \\
  \hline
40  & 0.935354 & 50 & 0.944776 & 60 & 0.951569 & 70  & 0.956728 & 100 & 0.966834   \\
  \hline
\end{tabular}
\end{center}
\caption{$R_k$ is the unique root of the equation
$  4(1-r)-r^{k-1}(1-2r+5r^2)=0 $  in $(0,1)$\label{Table2}}
\end{table}

\bthm\label{PW1-th2}
Let $f\in \mathcal{B}_k, k\geq 2,$   have an expansion \eqref{a}. Then the inequality
\[\sum_{n=k}^{\infty}|a_n|r^n+\left(\frac{r^{-k}}{1+|a_k|}+\frac{r^{1-k}}{1-r}\right)\sum_{n=k}^{\infty}|a_n|^2r^{2n}\leq 1
\]
is valid for $r\in [0,S_k]$, where $S_k$ is the unique solution in  $(0,1)$ of the equation
\[  2(1-r)-r^k(3-r)=0.
\]
The upper bound $S_k$ cannot be improved.
\ethm

From Theorems \ref{PW1-th1} and \ref{PW1-th2}, it is clear that $S_k\leq R_k$.
Numerical values of $S_k$ for certain values of $k\geq 2$ are listed in Table \ref{Table3}. Also, it is worth pointing out
that $\rho_k(1) = S_k$, where  $\rho_k(a)$ as in Theorem \ref{PW1-th3}.

\begin{table}
\begin{center}
\begin{tabular}{|l|l||l|l||l|l||l|l||l|l||}
\hline
 $k$   & $S_k$ & $k$ &  $S_k$  & $k$ & $S_k$ & $k$   & $S_k$ & $k$ &  $S_k$\\
\hline
2  & 0.585786 & 3 & 0.66152 & 4 &  0.709616 & 5  & 0.743563 & 6 & 0.769115  \\
  \hline
7  & 0.789207 & 8 & 0.805514 & 9 &  0.819071 & 10  & 0.830558 & 11 & 0.840442   \\
  \hline
15  & 0.869417 & 20 & 0.892242 & 25 & 0.907521 & 30  & 0.918574 & 35 & 0.926998   \\
  \hline
40  & 0.933662 & 50 & 0.943594 & 60 & 0.950691 & 70  & 0.956047 & 100 & 0.966459  \\
  \hline

\end{tabular}
\end{center}
\caption{$S_k$ is the unique root of the equation
$2(1-r)-r^k(3-r)=0$  in $(0,1)$ \label{Table3}}
\end{table}

\bthm\label{PW1-th3}
Let $f\in \mathcal{B}_k, k\geq 2,$  have an expansion  \eqref{a},  and let $|a_k|=a\in (0,1]$  be fixed. Then the inequality
\[\sum_{n=k}^{\infty}|a_n|r^n+\left(\frac{r^{-k}}{1+|a_k|}+\frac{r^{1-k}}{1-r}\right)\sum_{n=k}^{\infty}|a_n|^2r^{2n}\leq 1
\]
is valid for $r\in [0,\rho_k(a)]$,   where  $\rho_k(a)$ is the unique  solution in  $(0,1)$   of the equation
\[  (1+a)(1-r)-r^k[2a^2+a+r(1-2a^2)]=0.
\]
 The upper bound $\rho_k(a)$ cannot be improved.
\ethm

We see that $\rho_k(a)$ is a decreasing function of $a$. Also, we find that $\rho_2(\frac{1}{\sqrt{2}})=\frac{\sqrt{5}-1}{2}.$

Numerical values of $\rho_k(a)$ for certain values of $k\geq 2$ for fixed $a\in (0,1]$ are listed in Tables \ref{Table4(a)}- \ref{Table4(d)}.
\begin{table}
\begin{center}
\begin{tabular}{|l|l||l|l||l|l||l|l||l|l||}
\hline
 $k$   & $\rho_k(5/6)$ & $k$ &  $\rho_k(5/6)$  & $k$ & $\rho_k(5/6)$ & $k$   & $\rho_k(5/6)$ & $k$ &  $\rho_k(5/6)$\\
  \hline
2  & 0.604242 & 3 & 0.673433 & 4 &  0.718134 & 5  & 0.750042 & 6 & 0.774255 \\
  \hline
7  & 0.79341 & 8 & 0.80903 & 9 &  0.822067 & 10  & 0.833149 & 11 & 0.842709   \\
  \hline
15  & 0.870869 & 20 & 0.89319 & 25 & 0.908196 & 30  & 0.919083 & 35 & 0.927398  \\
  \hline
40  & 0.933985 & 50 & 0.943819 & 60 & 0.950858 & 70  & 0.956177 & 100 & 0.966531  \\
  \hline

\end{tabular}
\end{center}
\caption{$\rho_k(5/6)$ is the unique root of the equation
$\frac{11 (1 - r)}{6} - (\frac{20}{9}   -\frac{7r}{8} ) r^k=0$  in $(0,1)$\label{Table4(a)}}
\end{table}

\begin{table}
\begin{center}
\begin{tabular}{|l|l||l|l||l|l||l|l||l|l||}
\hline
 $k$   & $\rho_k(3/4)$ & $k$ &  $\rho_k(3/4)$  & $k$ & $\rho_k(3/4)$ &  $k$   & $\rho_k(3/4)$ & $k$ &  $\rho_k(3/4)$\\
  \hline
2  & 0.613378 & 3 & 0.679324 & 4 &  0.722344 & 5  & 0.753244 & 6 & 0.776794 \\
  \hline
7  & 0.795485 & 8 & 0.810767 & 9 &  0.823546 & 10  & 0.834427 & 11 & 0.843827  \\
  \hline
15  & 0.871585 & 20 & 0.893657 & 25 & 0.908528 & 30  & 0.919334 & 35 & 0.927594 \\
  \hline
40  & 0.934144 & 50 & 0.94393 & 60 & 0.950941 & 70  & 0.956241 & 100 & 0.966566 \\
  \hline

\end{tabular}
\end{center}
\caption{$\rho_k(3/4)$ is the unique root of the equation
$\frac{7 (1 - r)}{4} - (\frac{15}{8}   -\frac{r}{8} ) r^k=0$  in $(0,1)$\label{Table4(b)}}
\end{table}

\begin{table}
\begin{center}
\begin{tabular}{|l|l||l|l||l|l||l|l||l|l||}
\hline
$k$   & $\rho_k(2/3)$ & $k$ &  $\rho_k(2/3)$  & $k$ & $\rho_k(2/3)$ & $k$   & $\rho_k(2/3)$ & $k$ &  $\rho_k(2/3)$\\
  \hline
2  & 0.622387 & 3 & 0.685138 & 4 &  0.726502 & 5  & 0.756407 & 6 & 0.779302\\
  \hline
7  & 0.797536 & 8 & 0.812481 & 9 &  0.825007 & 10  & 0.835689 & 11 & 0.844932 \\
  \hline
15  & 0.872292 & 20 & 0.894118 & 25 & 0.908856 & 30  & 0.919581 & 35 & 0.927788 \\
  \hline
40  & 0.934301 & 50 & 0.944039 & 60 & 0.951022 & 70  & 0.956304 & 100 & 0.966601 \\
  \hline

\end{tabular}
\end{center}
\caption{$\rho_k(2/3)$ is the unique root of the equation
$\frac{5 (1 - r)}{3} - (\frac{14}{9}   + \frac{r}{9} ) r^k=0 $  in $(0,1)$\label{Table4(c)}}
\end{table}

\begin{table}
\begin{center}
\begin{tabular}{|l|l||l|l||l|l||l|l||l|l||}
\hline
$k$   & $\rho_k(1/2)$ & $k$ &  $\rho_k(1/2)$  & $k$ & $\rho_k(1/2)$ & $k$   & $\rho_k(1/2)$ & $k$ &  $\rho_k(1/2)$\\
  \hline
2  & 0.639802 & 3 & 0.696418 & 4 &  0.734582 & 5  & 0.762559 & 6 & 0.784184\\
  \hline
7  & 0.801527 & 8 & 0.815821 & 9 &  0.827851 & 10  & 0.838148 & 11 & 0.847082 \\
  \hline
15  & 0.873668 & 20 & 0.895014 & 25 & 0.909493 & 30  & 0.92006 & 35 & 0.928164 \\
  \hline
40  & 0.934605 & 50 & 0.944251 & 60 & 0.951179 & 70  & 0.956426 & 100 & 0.966668 \\
  \hline

\end{tabular}
\end{center}
\caption{$\rho_k(1/2)$ is the unique root of the equation
$\frac{3(1 - r)}{2} - (1   + \frac{r}{2} ) r^k=0$  in $(0,1)$\label{Table4(d)}}
\end{table}

\section{Proof of the Theorems}
\subsection{Proof of the Theorem \ref{PW1-th1}}
 Since $f\in \mathcal{B}_k,$ according to the classical lemma of Schwarz we may write $f(z)=z^kg(z)$, where $g\in \mathcal{B}.$ Let
\[  g(z)=\sum_{n=0}^{\infty}b_nz^n.
\]
Obviously, we have $a_{n+k}=b_n$ for $n\geq 0$. If we apply Theorem \Ref{PVW1-theC} to the function $g$, insert the last equation,
and multiply this inequality by $r^k$, we get
\be\label{c}
\sum_{n=k}^{\infty}|a_n|r^n+\left(\frac{1}{1+|a_k|}+\frac{r}{1-r}\right)\sum_{n=k+1}^{\infty}|a_n|^2r^{2n-k}\leq r^k\left(|a_k|+\frac{r}{1-r}(1-|a_k|^2)\right).\ee
If we fix $r$ in the expression on the right hand side and calculate the maximum with respect to $|a_k|$, we arrive
at the maximum value $r^kM(r)$ which is achieved for $|a_k|=1$, if $r\in [0,1/3]$, and for $|a_k|=(1-r)/2r$ in the remaining cases.
This value is less than or equal to unity if and only if $s_k(r)\geq 0$, where
\[  s_k(r)=4(1-r)-r^{k-1}(1-2r+5r^2).
\]
This proves the first part of the assertion. To prove the uniqueness of the solution in $(0,1)$ of  $s_k(r)=0$,
we firstly consider the case $k=2$. We see that $s_2(0)>0$, $s_2(1)<0,$ and
\[  s_2'(r)=-(5-4r+15r^2)<0 ~\mbox{ for }~ r\in [0,1].
\]
Since $s_2$ is monotonically decreasing on $[0,1]$, the assertion is proved.

For $k\geq 3$, we have $s_k(0)>0$, $s_k(1)<0,$ and $s_k'(0)=-4$. Further,
\[  s_k''(r)=-r^{k-3}\left[(k-1)(k-2)-2k(k-1)r+5(k+1)kr^2\right ].
\]
The discriminant of the quadratic form in $r$ in the square bracket term is negative for all $k\geq 3$, and thus, $ s_k''(r)<0$ so that $ s_k'(r)\leq  s_k'(0)<0$.
Hence, both $s_k'$ and $s_k$ are monotonically decreasing on $(0,1)$. This proves the uniqueness for $k\geq 3$.

For the proof of the sharpness of the upper bound we use the function
\be\label{e}
f(z)=z^k\left(\frac{a-z}{1-az}\right)=az^k-(1-a^2)\sum_{n=1}^{\infty}a^{n-1}z^{k+n},\quad a\in [0,1].
\ee
In this case we get
\[  \sum_{n=k}^{\infty}|a_n|r^n=ar^k+(1-a^2)\frac{r^{k+1}}{1-ar}
\]
and
\[\left(\frac{1}{1+|a_k|}+\frac{r}{1-r}\right)\sum_{n=k+1}^{\infty}|a_n|^2r^{2n-k}=\frac{(1-a^2)(1-a)r^{k+2}}{(1-r)(1-ar)}.
\]
The sum of these two terms is equal to
\[  ar^k+\frac{r^{k+1}(1-a^2)}{1-r}
\]
which equals unity for $a=(1-R_k)/2R_k.$ The proof of Theorem \ref{PW1-th1} is completed. \hfill $\Box$

\subsection{Proof of Corollary \ref{PW1-cor1}} According to Theorem \ref{PW1-th1}, equality in \eqref{b} for $r=R_k$ can be achieved
if and only if $|a_n|=0$ for $n\geq k+1.$ If this is assumed, we get
\[ \sum_{n=k}^{\infty}|a_n|r^n=|a_k|r^k<1
\]
for $r=R_k$ and a forteriori for $r\in [0,R_k].$ \hfill $\Box$

\br
Let us consider the cases $k\geq 2$. In Theorem \Ref{PVW1-theB} there are two possibilities. Either
\[
  \frac{r_k^k(1-2r_k+5r_k^2)}{4r_k(1-r_k)} = 1,\]
which means $r_k=R_k$, or
\[\frac{r_k^k}{\sqrt{1-r_k^2}} = 1.
\]
If the radius $r_k$ in the first possibility would be sharp, there would exist a uniformly convergent sequence
\[ f_m(z) = \sum_{n=k}^{\infty}a_{m,n}z^n\,\in \mathcal{B}_k,
\]
such that
\begin{equation}\label{Pro}
  \lim_{m\to \infty}\sum_{n=k}^{\infty}|a_{m,n}|R_k^n = 1.
\end{equation}
Analogous to the proof of Corollary \ref{PW1-cor1} this leads to
\[\lim_{m\to\infty}|a_{m,n}| = 0, \quad  n\geq k+1,
\]
a contradiction to \eqref{Pro}. This rules out the sharpness in the first
possibility, which answers in part the question of Paulsen et al. \cite{PaulPopeSingh-02-10}.

If there exists a function in $\mathcal{B}_k$ satisfying the condition  $|a_n|= r_k^n$ for $n\geq k$, then the second possibility in Theorem \Ref{PVW1-theB}
is sharp. This is because Cauchy-Schwarz inequality yields that
\be\label{g}
\sum_{n=k}^{\infty}|a_n|r^n\leq \sqrt{\sum_{n=k}^{\infty}|a_n|^2}\sqrt{\sum_{n=k}^{\infty}r^{2n}}\leq \frac{r^k}{(1-r^2)^{\frac{1}{2}}}
\leq 1
~\mbox{ for $r\in [0,r_k]$},
\ee
and equality in the first inequality of this chain  gives
\[     |a_n|=cr^n ~\mbox{ for $n\geq k$ with $c>0$}
\]
while equality in the second inequality implies $c=1$.
\er

\subsection{Proof of Theorem \ref{PW1-th2}} The inequality \eqref{c} implies
\be\label{d}
\sum_{n=k}^{\infty}|a_n|r^n+\left(\frac{r^{-k}}{1+|a_k|}+\frac{r^{1-k}}{1-r}\right)\sum_{n=k}^{\infty}|a_n|^2r^{2n}
\leq r^k\left(|a_k|+\frac{r}{1-r}+\frac{|a_k|^2}{1+|a_k|}\right).
\ee
Since for fixed value of $r$ the right hand side of this inequality is a monotonic increasing function of $|a_k|\in [0,1]$,
it is evident that the left hand side is less than or equal to unity, if and only if
\[r^k\left(\frac{3-r}{2(1-r)}\right)\leq 1.
\]
The equation
\[   t_k(r):= 2(1-r)-r^k(3-r)=0
\]
has a unique solution in (0,1), since $t_k(0)>0$, $t_k(1)<0$ and $t_k'(r) < 0$ for $r\in [0,1].$
This proves the first part of the assertion.

Further, it is easy to see that for $f(z)=z^k$ both sides of \eqref{d} are equal. Since we found the maximum of the right hand side for $|a_k|=1$,
it is obvious that this function proves sharpness of the upper bound.   The proof of Theorem \ref{PW1-th2} is completed.
\hfill $\Box$

\subsection{Proof of Theorem \ref{PW1-th3}} From \eqref{d} it is obvious that the left hand side of \eqref{d} is less than or equal to unity if
 \be\label{f}
 r^k\left(a+\frac{r}{1-r}+\frac{a^2}{1+a}\right)=\frac{r^k\left [2a^2+a+r(1-2a^2)\right ]}{(1+a)(1-r)}\leq 1.
 \ee
 This proves the first part of the assertion.
 Next, we have to prove the uniqueness of the solution of
 \[ u_{k,a}(r):=(1+a)(1-r)-r^k\left [2a^2+a+r(1-2a^2)\right ]=0
 \]
 in $[0,1]$. We use a method analogous to the above ones. It is obvious that $u_{k,a}(0)>0$ and $u_{k,a}(1)<0$. Further,
 \[u_{k,a}'(r)=-(1+a)-r^{k-1}k(2a^2+a)-r^k(k+1)(1-2a^2),\]
 and
 \[u_{k,a}''(r)=-kr^{k-2}A(r), \quad A(r)= (k-1)(2a^2+a)-r(k+1)(2a^2-1) .
 \]
 For $a\in [0,1/\sqrt{2}]$, it is clear that $u_{k,a}'(r) < 0$ for $r\in [0,1].$ In the remaining cases, we see that the
 second derivative is non-positive using the estimate
 \[A(r) \geq A(1)= -4a^2+(k-1)a+(k+1) \geq -4a^2+a+3 =(1-a)(4a+3) \geq 0,
 \]
 for $a\in [1/\sqrt{2},1]$ and $k\geq 2$. This together with $u_{k,a}'(0)<0$ ensures the negativity of the first derivative and
 in turn the uniqueness of the solution $\rho_{k}(a)$.

 To prove the sharpness of this upper bound, we again use the function \eqref{e}, but here we fix $a\in [0,1]$. In this case we get
 for the left hand side of \eqref{d} (for simplicity call it as $B(r)$) the term
 \beqq
 B(r) &=&    ar^k+\frac{r^{k+1}(1-a^2)}{1-ar}+\frac{r^{-k}(1+ar)}{(1+a)(1-r)}\left(a^2r^{2k}+\frac{r^{2k+2}(1-a^2)^2}{1-a^2r^2}\right)\\
 &=& r^k\left [a+\frac{a^2(1+ar)}{(1+a)(1-r)}+\frac{r(1-a^2)}{1-ar}\left(1+\frac{r(1-a)}{1-r}\right)\right ]\\
 &= & \frac{r^k[2a^2+a+r(1-2a^2)]}{(1+a)(1-r)}.
 \eeqq
  Comparison of  this expression with the right hand side of the equation in formula \eqref{f} delivers the asserted sharpness.
  \hfill $\Box$

  \bigskip
  \noindent
  Added in proof: Theorem 1.1 in \cite{BomBor-04} shows that the possibility discussed in Remark 2 cannot occur. Hence, the answer to the above mentioned question from \cite{PaulPopeSingh-02-10} is in fact negative.



\begin{thebibliography}{99}

\bibitem{AAP2016} Y. Abu Muhanna, R. M. Ali, and S. Ponnusamy,
On the Bohr inequality,
\emph{In "Progress in Approximation Theory and Applicable Complex Analysis" (Edited by N.K. Govil et al. ),
Springer Optimization and Its Applications} \textbf{117} (2016), 265--295.

\bibitem{AlKayPON-19} S. A. Alkhaleefah, I. R. Kayumov and S. Ponnusamy,
On the Bohr inequality with a fixed zero coefficient,
\emph{Proc. Amer. Math. Soc.} \textbf{147}(12) (2019), 5263--5274.

\bibitem{BDK5} C. B\'en\'eteau, A. Dahlner and D. Khavinson,
Remarks on the Bohr phenomenon,
\emph{Comput. Methods Funct. Theory} {\bf 4}(1) (2004),  1--19.

\bibitem{BhowDas-18} B. Bhowmik and N. Das,
Bohr phenomenon for subordinating families of certain univalent functions,
\emph{J. Math. Anal. Appl.} \textbf{462}(2) (2018), 1087--1098.

\bibitem{BoasKhavin-97-4} H. P. Boas and  D. Khavinson,
Bohr's power series theorem in several variables,
\emph{Proc. Amer. Math. Soc.} \textbf{125}(10) (1997),  2975--2979.

\bibitem{Bohr-14} H. Bohr,
A theorem concerning power series,
\emph{Proc. London Math. Soc.} \textbf{13}(2) (1914), 1--5.

\bibitem{Bom-62}
 E. Bombieri, Sopra un teorema di H. Bohr e G. Ricci sulle funzioni maggioranti delle serie
di potenze,
\emph{Boll. Un. Mat. Ital.} \textbf{17} (3)(1962), 276--282.

\bibitem{BomBor-04} E. Bombieri and J. Bourgain, A remark on Bohr's inequality,
\emph{Int. Math. Res. Not.} \textbf{80} (2004), 4307--4330.


\bibitem{Carlson-40}
F.  Carlson,  Sur les coefficients d'une fonction born\'{e}e dans le cercle unit\'{e} (French)
\emph{Ark. Mat. Astr. Fys.} \textbf{27A}(1) (1940), 8 pp.

\bibitem{EvPoRa-19} S. Evdoridis, S. Ponnusamy and A. Rasila,
Improved Bohr's inequality for locally univalent harmonic mappings,
\emph{Indag. Math. (N.S.)}  \textbf{30} (2019),  201--213.

\bibitem{EvPoRa-20} S. Evdoridis, S. Ponnusamy and A. Rasila,
Improved Bohr's inequality for mappings defined on simply connected domains, (in preparation).


\bibitem{FR} R. Fournier and S. Ruscheweyh,
On the Bohr radius for simply connected plane domains. {\it Hilbert spaces of analytic functions}, 165--171,
CRM Proc. Lecture Notes, 51, Amer. Math. Soc., Providence, RI, 2010.

\bibitem{GarMasRoss-2018}
S.~R.~Garcia, J.~ Mashreghi and W.~T.~Ross,
\emph{Finite Blaschke products and their connections}, Springer, Cham, 2018.


\bibitem{KayPon1} I. R. Kayumov and S. Ponnusamy,
Bohr inequality for odd analytic functions,
\emph{Comput. Methods Funct. Theory} \textbf{17} (2017), 679--688.

\bibitem{KayPon3} I. R. Kayumov and S. Ponnusamy, Improved version of Bohr's inequality,
\emph{Comptes Rendus Mathematique}
\textbf{356}(3) (2018),  272--277


\bibitem{KayPon2} I. R. Kayumov and S. Ponnusamy,
Bohr's inequalities for the analytic functions with lacunary series and harmonic functions,
\emph{J. Math. Anal. and Appl.,}  \textbf{465} (2018), 857--871.

\bibitem{KP-AASFM2-19} I. R. Kayumov  and S. Ponnusamy, On a powered Bohr inequality,
\emph{Ann. Acad. Sci. Fenn.  Ser. A I Math.}  \textbf{44}(2019), 301--310.


\bibitem{Landau-29} 
E. Landau,
Darstellung und Begr\"{u}ndung einiger neuerer Ergebnisse der Funktionentheorie,
Springer-Verlag, 1986

\bibitem{PaulPopeSingh-02-10} V. I. Paulsen, G. Popescu and D. Singh,
On Bohr's inequality,
\emph{Proc. London Math. Soc. } {\bf 85}(2) (2002), 493--512.

\bibitem{PauSin--09}
V. I. Paulsen and D. Singh, Bohr's inequality for uniform algebras,
\emph{Proc. Amer. Math. Soc.} \textbf{132(12)}(2004), 3577--3579,

\bibitem{PVW} S. Ponnusamy, R. Vijayakumar, and K.-J. Wirths, Refinement of the classical Bohr inequality, \emph{Preprint.}\\
{\tt http://arxiv.org/abs/1911.05315}
 .

%

\bibitem{Ricci-1956} G. Ricci, Complementi a un teorema di H. Bohr riguardante le serie di potenze,
\emph{Rev. Un.Mat. Argentina} \textbf{17} (1955/1956), 185--195.

%
%
\bibitem{Rogo-43}  W. Rogosinski, On the coefficients of subordinate functions,
\emph{Proc. London Math. Soc.} \textbf{48}(2) (1943), 48--82.

\bibitem{Sidon-27-15} S. Sidon, \"{U}ber einen Satz von Herrn Bohr,
\emph{Math. Z.} \textbf{26}(1) (1927), 731--732.

\bibitem{Tomic-62-16} M. Tomi\'c, Sur un th\'eor\`eme de H. Bohr,
\emph{Math. Scand.} \textbf{11} (1962), 103--106.

\end{thebibliography}
\end{document}